\documentclass[a4paper,11pt]{article}
\usepackage{bm,mathrsfs,amsmath,amsthm,amssymb,amscd,longtable,array,graphicx}
\newcommand{\nc}{\newcommand}
\nc{\rnc}{\renewcommand}
\nc{\nn}{\nonumber}
\nc{\der}{{\partial}}
\rnc{\Im}{{\rm{Im}\,}}
\rnc{\Re}{{\rm{Re}\,}}
\nc{\db}{\displaybreak[0]\\}
\nc{\bra}{\langle}
\nc{\ket}{\rangle}
\nc{\bs}{\boldsymbol}

\newtheorem{theorem}{Theorem}[section]

\newtheorem{proposition}[theorem]{Proposition}

\theoremstyle{definition}

\numberwithin{equation}{section}

\numberwithin{equation}{section}

\begin{document}%
%
\title{Integrability approach to
Feh\'er-N\'emethi-Rim\'anyi-Guo-Sun type identities
for factorial Grothendieck polynomials}

\author{
Kohei Motegi \thanks{E-mail: kmoteg0@kaiyodai.ac.jp}
\\\\
{\it Faculty of Marine Technology, Tokyo University of Marine Science and Technology,}\\
 {\it Etchujima 2-1-6, Koto-Ku, Tokyo, 135-8533, Japan} \\
\\\\
\\
}

\date{\today}

\maketitle

\begin{abstract}
Recently, Guo and Sun derived an identity for factorial
Grothendieck polynomials which is a generalization of the one
for Schur polynomials by Feh\'er, N\'emethi and Rim\'anyi.
We analyze the identity from the point of view of quantum integrability,
based on the correspondence between the wavefunctions of a five-vertex model
and the Grothendieck polynomials.
We give another proof using the quantum inverse scattering method.
We also apply the same idea and technique to derive an identity for
factorial Grothendieck polynomials for rectangular Young diagrams.
Combining with the Guo-Sun identity, we get a duality formula.
We also discuss a $q$-deformation of the Guo-Sun identity.
\end{abstract}

\section{Introduction}
Recently, Guo and Sun derived identities for
factorial Grothendieck polynomials \cite{GS}.
The factorial $\beta$-Grothendieck polynomials,
which is a $K$-theoretic analogue of the factorial Schur polynomials
\cite{LS,FK,Buch,Mc}, have the following determinant form
\cite{IN,IS}
\begin{align}
G_\lambda(\bs{z}|\bs{\alpha})=
   \frac{\mathrm{det}_n([z_i|\bs{\alpha}]^{\lambda_j+n-j}(1+\beta z_i)^{j-1})}
        {\prod_{1 \le i < j \le n}(z_i-z_j)},
 \label{factorialGrothendieck}
\end{align}
where $\lambda=(\lambda_1,\lambda_2,\dots,\lambda_n)$
is a  partition, i.e. a nonincreasing sequence of nonnegative integers
whose graphical representation is naturally given by the Young diagram.
$\bs{z}=\{z_1, \dots, z_n \}$ is a set of symmetric variables,
$\bs{\alpha}=\{\alpha_1,\alpha_2,\dots \}$ is a set of factorial variables
and
\begin{align}
[z_i|\bs{\alpha}]^j=(z_i \oplus \alpha_1) (z_i \oplus \alpha_2)
\cdots (z_i \oplus \alpha_j),
\end{align}
where $z \oplus \alpha:=z+\alpha+\beta z \alpha$.

One of the identities Guo and Sun derived
is the following one \cite{GS}:
for a partition $\lambda=(\lambda,\dots,\lambda_k)$ such that
$\lambda_1 \le m-k$ and another one $\mu=((m-k)^{n-k},\lambda_1,\dots,\lambda_k)$,
the following identity holds:
\begin{align}
G_{\mu}(\bs{z}|\bs{\alpha})=\sum_{S_k^n \in \binom{[n]}{k}}
G_\lambda(\bs{z}_{S_k^n}|\bs{\alpha})
\frac{\displaystyle \prod_{i \in {S_k^n}} (1+\beta z_i)^{n-k} \prod_{j \in \overline{S_k^n}} [z_j|\bs{\alpha}]^m}{\displaystyle \prod_{i \in {S_k^n}} \prod_{j \in \overline{{S_k^n}}}(z_j-z_i)}, \label{guosun}
\end{align}
where $S_k^n$ is a $k$-subset of $[n]=\{1,2,\dots,n \}$,
$\binom{[n]}{k}$ is the set of $k$-subsets of $n$,
$\overline{S_k^n}=\{1,2,\dots,n \} \backslash S_k^n$,
and
$G_\lambda(\bs{z}_{S_k^n}|\bs{\alpha})=G_\lambda(\{ z_{i_1},\dots,z_{i_k} \}|\bs{\alpha})$ for $S_k^n=\{i_1,\dots,i_k \}$.

This identity generalizes the one for the Schur polynomials
derived by Feh\'er, N\'emethi and Rim\'anyi \cite{FNR}
which corresponds to the case
$\beta=\alpha_1=\alpha_2=\dots=0$,
and we will call this type of identity as
Feh\'er-N\'emethi-Rim\'anyi-Guo-Sun type identity.

In this paper, we restrict to the case $\beta=-1$,
and investigate this identity and also derive similar identities
from the viewpoint of quantum integrability.
Recently, there is an active line of research which investigates relations
between integrable models and related structures (integrable lattice models,
classical integrable systems, vertex operators, crystal basis)
and the (dual, symmetric) Grothendieck polynomials,
and study the properties of the
Grothendieck polynomials using the connections.
See \cite{MS,MS2,KirillovSigma,GK,WZ,IwNa,Iwao,MOS,BSW}
for examples for various topics.
We give another proof of the Guo-Sun identity \eqref{guosun}
using the quantum inverse scattering method \cite{FST,KBI},
which is a method developed to study quantum integrable models.
Why we can use this method for giving another proof
is based on the correspondence between certain types of partition functions
of an integrable five-vertex model and the Grothendieck polynomials.
This correspondence was used to investigate Cauchy-type identities,
Gromov-Witten invariants and Littlewood-Richardson coefficients
\cite{MS,MS2,GK,WZ}.
In this paper, we give another application of the correspondence.
Namely, we give an integrability proof of the Guo-Sun identity
using the quantum inverse scattering method.
We also apply the same idea and technique to derive
an identity for factorial Grothendieck polynomials
for rectangular Young diagrams.
The five-vertex model which is used can be regarded as a certain limit
of the $U_q(\widehat{sl_2})$ six-vertex model \cite{Dr,J,LW}.
Based on this viewpoint, we also discuss a $q$-deformation
by following the line of computation to derive the Guo-Sun identity.

This paper is organized as follows.
In the next section, we explain the
correspondence between the wavefunctions of the
$U_q(\widehat{sl_2})$ six-vertex model and symmetric functions,
and its $q=0$ degeneration which gives
the correspondence between the wavefunctions of the five-vertex model
and the factorial Grothendieck polynomials.
In section 3, we give another proof 
of the Guo-Sun identity by using the quantum inverse scattering method.
In section 4, we derive an identity for the factorial Grothendieck polynomials for rectangular shapes by using the same idea and technique in section 3
to ``a different direction".
We discuss a $q$-deformation of the Guo-Sun identity by applying
the same idea to the
$U_q(\widehat{sl_2})$ six-vertex model.
Section 5 is devoted to the conclusion of this paper.

\section{$U_q(\widehat{sl_2})$ six-vertex model,
wavefunctions and degeneration to the five-vertex model
}

In this section, we first introduce the $U_q(\widehat{sl_2})$ six-vertex model
and explain the correspondence between the wavefunctions
and symmetric functions.
Next, we degenerate the six-vertex model to the five-vertex model
and explain the correspondence between the wavefunctions
and the factorial Grothendieck polynomials.

\begin{figure}[ht]
\includegraphics[width=12cm]{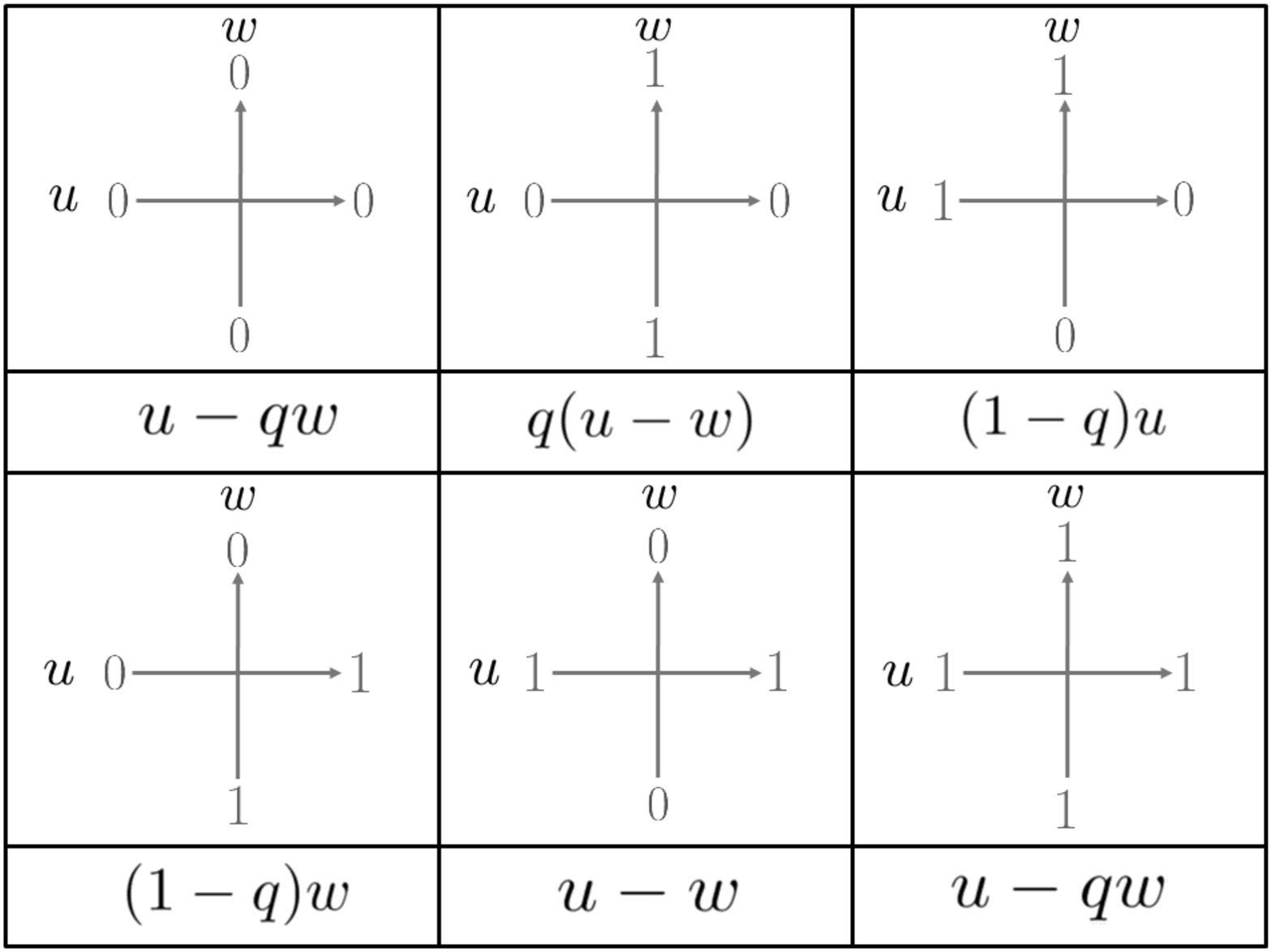}
\caption{The $R$-matrix  of the six-vertex model \eqref{rmatrix}.
There are six admissible vertex configurations,
and the weight associated to each configuration is given below
the figure of the configuration.}
\label{picturermatrix}
\end{figure}

The $U_q(\widehat{sl_2})$ $R$-matrix is the following matrix \cite{Dr,J}
(Figure \ref{picturermatrix})
\begin{eqnarray}
R_{ab}(u,w)=\left( 
\begin{array}{cccc}
u-qw & 0 & 0 & 0 \\
0 & q(u-w) & (1-q)u & 0 \\
0 & (1-q)w & u-w & 0 \\
0 & 0 & 0 & u-qw
\end{array}
\right), \label{rmatrix}
\end{eqnarray}
acting on the tensor product $W_a \otimes W_b$
of the complex two-dimensional space $W_a$.
We denote the dual space of $W_a$ by $W_a^*$.

The $R$-matrix \eqref{rmatrix} satisifies the Yang-Baxter relation
\begin{align}
R_{ab}(u,v)R_{ac}(u,w)R_{bc}(v,w)
=R_{bc}(v,w)R_{ac}(u,w)R_{ab}(u,v), \label{yangbaxter}
\end{align}
acting on $W_a \otimes W_b \otimes W_c$.

We denote the orthonormal basis of $W_a$ and its dual as
$\{|0 \rangle_a, |1 \rangle_a \}$ and $\{{}_a \langle 0|, {}_a \langle 1|\}$.
We also introduce
the following Pauli spin operators
$\sigma^+$ and $\sigma^-$ as operators acting on the (dual) orthonomal
basis as
\begin{align}
&\sigma^+|1 \rangle=|0 \rangle, \ 
\sigma^+|0 \rangle=0, \ 
\langle 0|\sigma^+=\langle 1|, \
\langle 1|\sigma^+=0, 
\\
&\sigma^-|0 \rangle=|1 \rangle, \
\sigma^-|1 \rangle=0, \
\langle 1|\sigma^-=\langle 0|, \
\langle 0|\sigma^-=0.
\end{align}

The monodromy matrix $T_a(u|w_1,\dots,w_{m+n-k})$ is the product of $R$-matrices
\begin{align}
T_{a}(u|w_1,\dots,w_{m+n-k})&=R_{a, m+n-k}(u,w_{m+n-k}) \cdots R_{a 1}(u,w_1)
\nonumber \\
&=
\begin{pmatrix}
A(u|w_1,\dots,w_{m+n-k}) & B(u|w_1,\dots,w_{m+n-k})  \\
C(u|w_1,\dots,w_{m+n-k}) & D(u|w_1,\dots,w_{m+n-k})
\end{pmatrix}_{a},
\label{monodromy}
\end{align}
acting on
$W_a \otimes W_1 \otimes \cdots \otimes W_{m+n-k}$.

The $B$-operator is a matrix element of the monodromy matrix
$T_{a}(u|w_1,\dots,w_{m+n-k})$
\begin{align}
B(u|w_1,\dots,w_{m+n-k})= {}_a
\langle 0 |T_a(u|w_1,\dots,w_{m+n-k})| 1 \rangle_a,
\end{align}
acting on $W_1 \otimes \cdots \otimes W_{m+n-k}$.

Let us define the (dual) vacuum state as
$|\Omega \rangle_{m+n-k}:=|0 \rangle_1 \otimes \cdots \otimes
|0 \rangle_{m+n-k}
\in W_1 \otimes \cdots \otimes W_{m+n-k}$
(
${}_{m+n-k} \langle \Omega |:={}_1 \langle 0 | \otimes \cdots \otimes
{}_{m+n-k} \langle 0|
\in W_1^* \otimes \cdots \otimes W_{m+n-k}^*$)
and configuration vectors as
\begin{align}
_{m+n-k} \langle x_1 \cdots x_n|
&={}_{m+n-k} \langle \Omega |
\prod_{j=1}^n \sigma^+_{x_j}
\in W_1^* \otimes \cdots \otimes W_{m+n-k}^*
, \label{ordinarydualparticleconfiguration}
\end{align}
for
$1 \le x_1 < x_2 < \cdots < x_n \le m+n-k$.

\begin{figure}[ht]
\includegraphics[width=12cm]{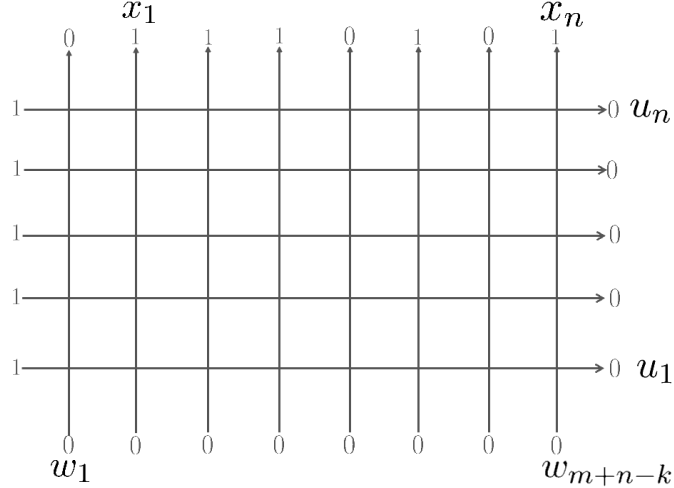}
\caption{The wavefunctions \eqref{sixvertexwavefunction}.
The sequence of 0s on the bottom part represents the vacuum state $|\Omega \rangle_{m+n-k}$. The $j$-th row counted from the bottom represents the $B$-operator
$B(u_j|w_1,\dots,w_{m+n-k})$ and the top part
which is a sequence of 0s and 1s represents
the dual vector $_{m+n-k} \langle x_1 \cdots x_n|$.
The dual vector is labelled by the positions $x_1,\dots,x_n$ of 1s counted from the left. The figure illustrates the case $n=5$, $m=k+3$,
$(x_1,x_2,x_3,x_4,x_5)=(2,3,4,6,8)$.
}
\label{picturewavefunction}
\end{figure}

We now introduce the wavefunctions
$W_{m+n-k,n}(u_1,\dots,u_n|w_1,\dots,w_{m+n-k}|x_1,\dots,x_n)$
as (Figure \ref{picturewavefunction})

\begin{align}
&W_{m+n-k,n}(u_1,\dots,u_n|w_1,\dots,w_{m+n-k}|x_1,\dots,x_n) \nonumber \\
=&_{m+n-k} \langle x_1 \cdots x_n|
B(u_n|w_1,\dots,w_{m+n-k}) \cdots B(u_1|w_1,\dots,w_{m+n-k})|\Omega \rangle_{m+n-k}.
\label{sixvertexwavefunction}
\end{align}

We next define the following symmetric function \\
$F_{m+n-k,n}(u_1,\dots,u_n|w_1,\dots,w_{m+n-k}|x_1,\dots,x_n)$
which depends on the symmetric variables $u_1,\dots,u_n$,
complex parameters $w_1,\dots,w_{m+n-k}$
and integers $x_1,\dots,x_n$ satisfying
$1 \le x_1 < \cdots < x_n \le m+n-k$,
\begin{align}
&F_{m+n-k,n}(u_1,\dots,u_n|w_1,\dots,w_{m+n-k}|x_1,\dots,x_n) \nonumber \\
=
&
\sum_{\sigma \in S_n}
\prod_{j=1}^n \prod_{i=x_j+1}^{m+n-k} (u_{\sigma(j)}-qw_i)
\prod_{1 \le i < j \le n}
\frac{qu_{\sigma(i)}-u_{\sigma(j)}}{u_{\sigma(i)}-u_{\sigma(j)}}
\nonumber \\
&\times
\prod_{j=1}^n \prod_{i=1}^{x_j-1}(u_{\sigma(j)}-w_i)
\prod_{j=1}^n (1-q)u_{\sigma(j)}.
\label{symmetricfunction}
\end{align}

The wavefunction
$W_{m+n-k,n}(u_1,\dots,u_n|w_1,\dots,w_{m+n-k}|x_1,\dots,x_n)$
is explicitly expressed as the
symmetric function
$F_{m+n-k,n}(u_1,\dots,u_n|w_1,\dots,w_{m+n-k}|x_1,\dots,x_n)$
\begin{align}
&W_{m+n-k,n}(u_1,\dots,u_n|w_1,\dots,w_{m+n-k}|x_1,\dots,x_n)
\nonumber \\
=&
F_{m+n-k,n}(u_1,\dots,u_n|w_1,\dots,w_{m+n-k}|x_1,\dots,x_n).
\label{correspondencesixvertex}
\end{align}
See \cite{BP,Motegi} for example for proofs of this correspondence.
Next, we explain the degeneration of the
correspondence \eqref{correspondencesixvertex}.
If one sets $q$ to $q=0$, the $R$-matrix
for the six-vertex model \eqref{rmatrix} reduces to
that for the five-vertex model
\begin{eqnarray}
R_{ab}(u,w)|_{q=0}=\left( 
\begin{array}{cccc}
u & 0 & 0 & 0 \\
0 & 0 & u & 0 \\
0 & w & u-w & 0 \\
0 & 0 & 0 & u
\end{array}
\right). \label{fivevertexrmatrix}
\end{eqnarray}
Under the change of variables
$z_j=1-u_j^{-1} \ (j=1,\dots,n)$, $\alpha_j=1-w_j \ (j=1,\dots, m+n-k)$,
$\lambda_j=x_{n-j+1}-n+j-1 \ (j=1,\dots,n)$,
the correspondence at $q=0$
\eqref{correspondencesixvertex}
becomes the
following correspondence between the
wavefunctions of the five-vertex model
and the $\beta=-1$ factorial Grothendieck polynomials \cite{MS,MS2,GK,WZ}
\begin{align}
&W_{m+n-k,n}(u_1,\dots,u_n|w_1,\dots,w_{m+n-k}|x_1,\dots,x_n)|_{q=0}
\nonumber \\
=&\prod_{j=1}^n \frac{1}{(1-z_j)^{m+n-k}}
G_\lambda(\{z_1,\dots,z_n \}|\bs{\alpha}).
\label{correspondence}
\end{align}
We use this correspondence in the next two sections
to investigate identities for the factorial Grothendieck polynomials
from the viewpoint of quantum integrability.

\section{Integrability proof of Guo-Sun identity}
In this and the next sections, we consider the five-vertex model
whose $R$-matrix is given by \eqref{fivevertexrmatrix}
which is the $q=0$ limit of \eqref{rmatrix}.
Every object introduced in the last section
should be understood that we set $q$ to
$q=0$ in this and the next sections.

In this section,
let us show another proof of Guo-Sun identity \eqref{guosun}
using the quantum inverse scattering method.
They derived an identity 
for Grothendieck polynomials
of the following partition $\mu=((m-k)^{n-k},\lambda_1,\dots,\lambda_k)$.
Applying the correspondence \eqref{correspondence}
to this partition,
we have
\begin{align}
&W_{m+n-k,n}(u_1,\dots,u_n|w_1,\dots,w_{m+n-k}|\lambda_k+1,\dots,\lambda_1+k,m+1,\dots,m+n-k)
\nonumber \\
=&\prod_{j=1}^n \frac{1}{(1-z_j)^{m+n-k}}
G_\mu(\{z_1,\dots,z_n \}|\bs{\alpha}).
\label{forcomparison}
\end{align}

\begin{figure}[ht]
\includegraphics[width=12cm]{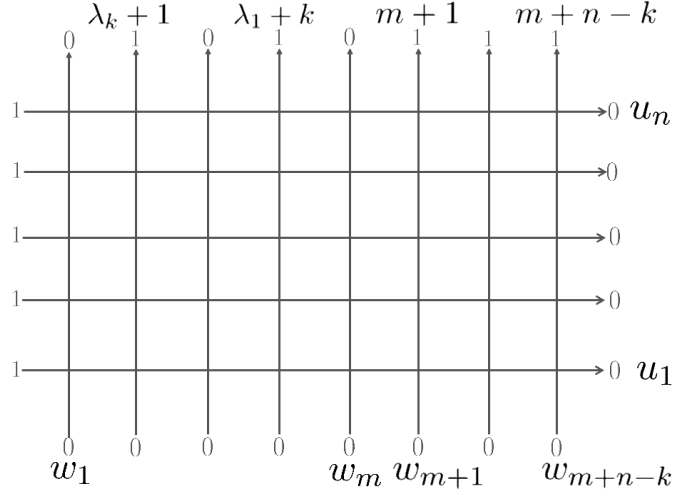}
\caption{The wavefunctions corresponding to the
factorial Grothendieck polynomials of partition
$\mu=((m-k)^{n-k},\lambda_1,\dots,\lambda_k)$ \eqref{forcomparison}.
The top part is a sequence of 0s and 1s.
The 1s appear at positions
$\lambda_k+1$, $\lambda_{k-1}+2$, $\dots$, $\lambda_1+k$, $m+1$, $m+2$, $\dots$, $m+n-k$ counted from the left. Subtracting $1,2,\dots,n$ from
$\lambda_k+1$, $\lambda_{k-1}+2$, $\dots$, $\lambda_1+k$, $m+1$, $m+2$, $\dots$, $m+n-k$, and reorder as a nonincreasing sequence,
we get the partition $((m-k)^{n-k},\lambda_1,\dots,\lambda_k)$
which labels the sequence of 0s and 1s in this figure.
The figure illustrates the case $n=5$, $m=5$, $k=2$,
$(\lambda_1,\lambda_2)=(2,1)$.
}
\label{picturewavefunctionguosun}
\end{figure}

To prove the Guo-Sun identity,
we investigate
$W_{m+n-k,n}(u_1,\dots,u_n|w_1,\dots,w_{m+n-k}|\lambda_k+1,\dots,\lambda_1+k,m+1,\dots,m+n-k)$ using its graphical description (Figure \ref{picturewavefunctionguosun}) and derive another expression.

\begin{figure}[ht]
\includegraphics[width=12cm]{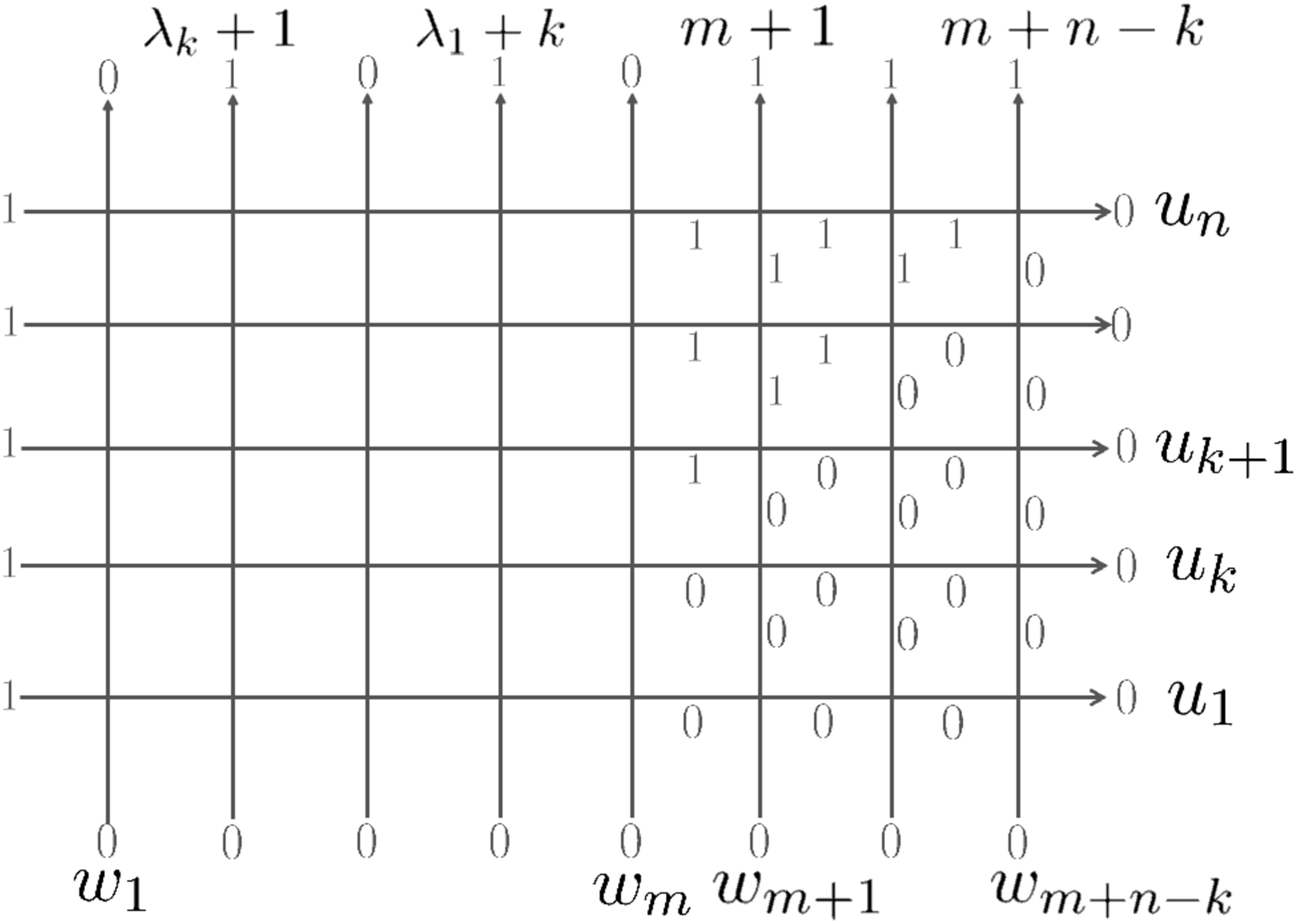}
\caption{
The wavefunction
$W_{m+n-k,n}(u_1,\dots,u_n|w_1,\dots,w_{m+n-k}|\lambda_k+1,\dots,\lambda_1+k,m+1,\dots,m+n-k)$.
We can see that
only one configuration is allowed in the rightmost $n-k$ columns,
which gives the factor $\displaystyle \prod_{j=1}^n u_j^{n-k}$.
The remaining part can be expressed as
$\displaystyle {}_{m} \langle \lambda_k+1,\dots,\lambda_1+k|
\prod_{j=k+1}^{n} D(u_j|w_1,\dots,w_m) \prod_{j=1}^k B(u_j|w_1,\dots,w_m)
|\Omega \rangle_m$.
}
\label{picturewavefunctionguosunfreeze}
\end{figure}

First, from the graphical representation of
$W_{m+n-k,n}(u_1,\dots,u_n|w_1,\dots,w_{m+n-k}|\lambda_k+1,\dots,\lambda_1+k,m+1,\dots,m+n-k)$ and noting that the bulk weights are given by the
$R$-matrix of the five-vertex model \eqref{fivevertexrmatrix},
one can see that
only one configuration is allowed in the rightmost $n-k$ columns
(Figure \ref{picturewavefunctionguosunfreeze}),
and we get the factor $\displaystyle \prod_{j=1}^n u_j^{n-k}$
from that configuration.
Next, looking at the remaining part, we find that they are expressed
as $\displaystyle {}_{m} \langle \lambda_k+1,\dots,\lambda_1+k|
\prod_{j=k+1}^{n} D(u_j|w_1,\dots,w_m) \prod_{j=1}^k B(u_j|w_1,\dots,w_m)
|\Omega \rangle_m$.
Hence we have
\begin{align}
&W_{m+n-k,n}(u_1,\dots,u_n|w_1,\dots,w_{m+n-k}|\lambda_k+1,\dots,\lambda_1+k,m+1,\dots,m+n-k)
\nonumber \\
=&\prod_{j=1}^n u_j^{n-k} {}_{m} \langle \lambda_k+1,\dots,\lambda_1+k|
\prod_{j=k+1}^{n} D(u_j|w_1,\dots,w_m) \prod_{j=1}^k B(u_j|w_1,\dots,w_m)
|\Omega \rangle_m. \label{tochu}
\end{align}
Next, we commute the multiple operators
$\prod_{j=k+1}^{n} D(u_j|w_1,\dots,w_m)$
with the multiple operators $\prod_{j=1}^k B(u_j|w_1,\dots,w_m)$.
From the Yang-Baxter relation
\eqref{yangbaxter}, we get the intertwining relation
for the monodromy matrices
\begin{align}
&R_{ab}(u_1,u_2)T_{a}(u_1|w_1,\dots,w_m)T_{b}(u_2|w_1,\dots,w_m)
\nonumber \\
=&T_{b}(u_2|w_1,\dots,w_m)T_{a}(u_1|w_1,\dots,w_m)R_{ab}(u_1,u_2).
\label{rttone}
\end{align}
Some matrix elements of \eqref{rttone} are given by
\begin{align}
D(u_1)B(u_2)&=\frac{u_1}{u_1-u_2}B(u_2)D(u_1)
-\frac{u_2}{u_1-u_2}B(u_1)D(u_2), \label{rttone1} \\
D(u_1)B(u_2)&=D(u_2)B(u_1), \label{rttone2} \\
B(u_1)B(u_2)&=B(u_2)B(u_1), \label{rttone3} \\
D(u_1)D(u_2)&=D(u_2)D(u_1). \label{rttone4}
\end{align}
For integrable models which have the $R$-matrix of the five-vertex model
as an intertwiner, a compact formula between the commutation relations
between multiple $D$-operators and $B$-operators
can be derived from \eqref{rttone1}, \eqref{rttone2}, \eqref{rttone3}, \eqref{rttone4}, following the argument of \cite{ShigechiUchiyama}
(they were analyzing the integrable phase model \cite{Bo},
but the type of the representation
for the quantum space does not affect the argument).
The result is given by
\begin{align}
&\prod_{j=k+1}^{n} D(u_j|w_1,\dots,w_m) \prod_{j=1}^k B(u_j|w_1,\dots,w_m)
\nonumber \\
=&\sum_{S_k^n \in \binom{[n]}{k}} \prod_{i \in S_k^n, j \in \overline{S_k^n}}
\frac{u_j}{u_j-u_i} \prod_{i \in S_k^n} B(u_i|w_1,\dots,w_m)
\prod_{j \in \overline{S_k^n}} D(u_j|w_1,\dots,w_m).
\label{multiplecommutation}
\end{align}
\eqref{multiplecommutation} can be shown as follows
(see the Proof of Theorem 6.1 in \cite{ShigechiUchiyama}).
First, using \eqref{rttone1}, \eqref{rttone3} and \eqref{rttone4}
to move all the $B$-operators to the left of all the $D$-operators,
one notes the operator part of all the terms which appear can be expressed as
\begin{align}
\prod_{i \in S_k^n} B(u_i|w_1,\dots,w_m)
\prod_{j \in \overline{S_k^n}} D(u_j|w_1,\dots,w_m),
\label{operatorpart}
\end{align}
for $\displaystyle S_k^n \in \binom{[n]}{k}$.
To extract the coefficient of \eqref{operatorpart}
for a fixed $S_k^n$,
one uses \eqref{rttone2} repeatedly to rewrite
the left hand side of \eqref{multiplecommutation} as
\begin{align}
\prod_{j \in \overline{S_k^n}} D(u_j|w_1,\dots,w_m)
\prod_{i \in S_k^n} B(u_i|w_1,\dots,w_m).
\end{align}
Finally, we only use \eqref{rttone1} repeatedly to move
all the $B$-operators to the left of all the $D$-operators.
We only need to concentrate on the first
term of the right hand side of \eqref{rttone1}
when commuting the $B$- and $D$-operators
to extract the coefficient
of \eqref{operatorpart},
since if we once use the second term of the left hand side of \eqref{rttone1},
we get other operators.
Noting this, one finds the coefficient of the operator is given by
$\displaystyle \prod_{i \in S_k^n, j \in \overline{S_k^n}} \frac{u_j}{u_j-u_i}$
, and we get the commutation relation \eqref{multiplecommutation}.

Using \eqref{multiplecommutation}
and the action of the $D$-operators on the vacuum state
\begin{align}
\prod_{j \in \overline{S_{k}^n}} D(u_j|w_1,\dots,w_m) | \Omega \rangle_m
=\prod_{j \in \overline{S_{k}^n}} \prod_{i=1}^m (u_j-w_i) | \Omega \rangle_m,
\end{align}
\eqref{tochu} becomes
\begin{align}
&W_{m+n-k,n}(u_1,\dots,u_n|w_1,\dots,w_{m+n-k}|\lambda_k+1,\dots,\lambda_1+k,m+1,\dots,m+n-k)
\nonumber \\
=&\prod_{j=1}^n u_j^{n-k}
\sum_{S \in \binom{[n]}{k}}
\prod_{i \in S_k^n, j \in \overline{S_k^n}}
\frac{u_j}{u_j-u_i}
\prod_{j \in \overline{S_{k}^n}} \prod_{i=1}^m (u_j-w_i)
\nonumber \\
\times&{}_{m} \langle \lambda_k+1,\dots,\lambda_1+k|
\prod_{i \in S_k^n} B(u_i|w_1,\dots,w_m)
|\Omega \rangle_m,
\end{align}
which one can further rewrite using
\eqref{sixvertexwavefunction} and \eqref{correspondence}
\begin{align}
{}_{m} \langle \lambda_k+1,\dots,\lambda_1+k|
\prod_{i \in S_k^n} B(u_i|w_1,\dots,w_m)
|\Omega \rangle_m
=\prod_{j \in S_k^n} \frac{1}{(1-z_j)^m}
G_\lambda(\bs{z}_{S_k^n}|\bs{\alpha}),
\end{align}
as
\begin{align}
&W_{m+n-k,n}(u_1,\dots,u_n|w_1,\dots,w_{m+n-k}|\lambda_k+1,\dots,\lambda_1+k,m+1,\dots,m+n-k)
\nonumber \\
=&\prod_{j=1}^n u_j^{n-k}
\sum_{S \in \binom{[n]}{k}}
\prod_{i \in S_k^n, j \in \overline{S_k^n}}
\frac{u_j}{u_j-u_i}
\prod_{j \in \overline{S_{k}^n}} \prod_{i=1}^m (u_j-w_i)
\prod_{j \in S_k^n} \frac{1}{(1-z_j)^m}
G_\lambda(\bs{z}_{S_k^n}|\bs{\alpha}). \label{tochu2}
\end{align}
Finally, we compare the two expressions for the same object.
From \eqref{forcomparison} and \eqref{tochu2},
we get
\begin{align}
&\prod_{j=1}^n \frac{1}{(1-z_j)^{m+n-k}}
G_\mu(\{z_1,\dots,z_n \}|\bs{\alpha})
\nonumber \\
=&\prod_{j=1}^n u_j^{n-k}
\sum_{S \in \binom{[n]}{k}}
\prod_{i \in S_k^n, j \in \overline{S_k^n}}
\frac{u_j}{u_j-u_i}
\prod_{j \in \overline{S_{k}^n}} \prod_{i=1}^m (u_j-w_i)
\prod_{j \in S_k^n} \frac{1}{(1-z_j)^m}
G_\lambda(\bs{z}_{S_k^n}|\bs{\alpha}),
\end{align}
which, using the translation rule
$z_j=1-u_j^{-1} \ (j=1,\dots,n)$, $\alpha_j=1-w_j \ (j=1,\dots, m+n-k)$,
can be rewritten as the Guo-Sun identity \eqref{guosun} for the case $\beta=-1$
\begin{align}
G_{\mu}(\bs{z}|\bs{\alpha})=\sum_{S_k^n \in \binom{[n]}{k}}
G_\lambda(\bs{z}_{S_k^n}|\bs{\alpha})
\frac{\displaystyle \prod_{i \in {S_k^n}} (1-z_i)^{n-k} \prod_{j \in \overline{S_k^n}} [z_j|\bs{\alpha}]^m}{\displaystyle \prod_{i \in {S_k^n}} \prod_{j \in \overline{{S_k^n}}}(z_j-z_i)}.
\end{align}

\section{An identity for rectangular shapes}
In this section, we apply the idea and technique
used in the last section ``in a different direction''
to derive an identity for factorial
Grothendieck polynomials.
We consider the case when the partitions whose corresponding Young diagrams are rectangular shapes,
i.e. we consider the case when the partition is of the form
$\mu=((m-k)^{n-k},0^k)$.
We show the following identity.

\begin{figure}[ht]
\includegraphics[width=12cm]{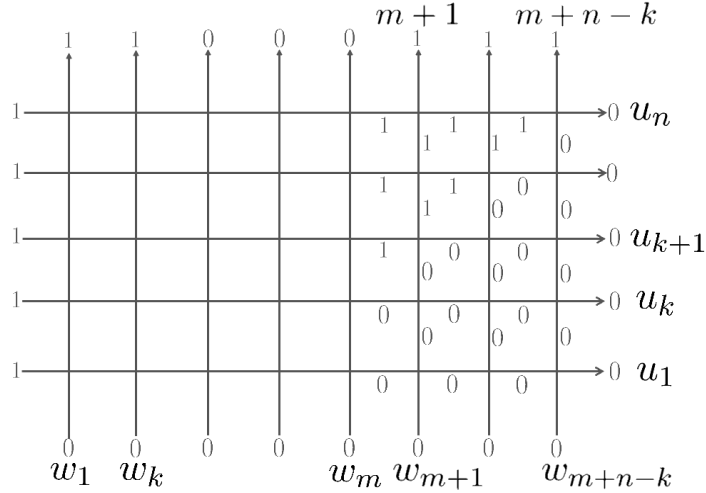}
\caption{The wavefunctions corresponding to the
factorial Grothendieck polynomials of 
rectangular shapes $\mu=((m-k)^{n-k},0^k)$
\eqref{forcomparisontwo}.
}
\label{picturewavefunctionrrectangular}
\end{figure}

\begin{theorem} \label{rectangulartheorem}
Let $\bs{z}=\{z_1, \dots, z_n \}$ be a set of symmetric variables
and
$\bs{\alpha}=\{\alpha_1,\alpha_2,\dots \}$ a set of factorial variables.
For a partition $\mu=((m-k)^{n-k},0^k)$,
the following identity holds:
\begin{align}
G_{\mu}(\bs{z}|\bs{\alpha})=\sum_{S_k^m \in \binom{[m]}{k}}
\frac{\displaystyle \prod_{i \in {S_k^m}} (1-\alpha_i)^{m-k} \prod_{j \in \overline{S_k^m}}
\prod_{i=1}^n (z_i \oplus \alpha_j)}
{\displaystyle \prod_{i \in {S_k^m}} \prod_{j \in \overline{{S_k^m}}}(\alpha_j-\alpha_i)}, \label{rectangularidentity}
\end{align}
where $S_k^m$ is a $k$-subset of $[m]=\{1,2,\dots,m \}$,
$\binom{[m]}{k}$ is the set of $k$-subsets of $m$ and
$\overline{S_k^m}=\{1,2,\dots,m \} \backslash S_k^m$.
\end{theorem}

\begin{proof}
We first introduce another class of monodromy matrices
\begin{align}
\overline{T}_j(w|u_1,\dots,u_n)&= R_{a_n j}(u_n|w) \cdots R_{a_1 j}(u_1|w)
\nonumber \\
&=
\begin{pmatrix}
\overline{A}(w|u_1,\dots,u_{n}) & \overline{B}(w|u_1,\dots,u_{n})  \\
\overline{C}(w|u_1,\dots,u_{n}) & \overline{D}(w|u_1,\dots,u_{n})
\end{pmatrix}_{j},
\label{anothermonodromy}
\end{align}
and vectors in the auxiliary space
$| \overline{\Omega} \rangle_n :=|1 \rangle_{a_1} \otimes \cdots \otimes
|1 \rangle_{a_n} \in W_{a_1} \otimes \cdots \otimes W_{a_n}$
and ${}_n \langle 0^k 1^{n-k}|:=
{}_{a_1} \langle 0 | \otimes \cdots \otimes {}_{a_k} \langle 0|
\otimes {}_{a_{k+1}} \langle 1 | \otimes \cdots \otimes {}_{a_n} \langle 1|
\in W_{a_1}^* \otimes \cdots \otimes W_{a_n}^*$.

In this section, we deal with the wavefunctions of the following type \\
$W_{m+n-k,n}(u_1,\dots,u_n|w_1,\dots,w_{m+n-k}|1,\dots,k,m+1,\dots,m+n-k)$,
which due to the correspondence \eqref{correspondence},
can be expressed using the factorial Grothendieck polynomials
of a rectangular shape $((m-k)^{n-k},0^k)$ as
\begin{align}
&W_{m+n-k,n}(u_1,\dots,u_n|w_1,\dots,w_{m+n-k}|1,\dots,k,m+1,\dots,m+n-k)
\nonumber \\
=&\prod_{j=1}^n \frac{1}{(1-z_j)^{m+n-k}}
G_{((m-k)^{n-k},0^k)}(\{z_1,\dots,z_n \}|\bs{\alpha}).
\label{forcomparisontwo}
\end{align}

We now apply the same idea and technique used in the last section, but we use
in a different direction this time. From the graphical representation
of  the wavefunctions $W_{m+n-k,n}(u_1,\dots,u_n|w_1,\dots,w_{m+n-k}|1,\dots,k,m+1,\dots,m+n-k)$
(Figure \ref{picturewavefunctionrrectangular}),
we first find that
only one configuration is allowed in the rightmost $n-k$ columns,
which gives the factor $\prod_{j=1}^n u_j^{n-k}$ again.
The remaining part can be written using the matrix elements of another
type of monodromy matrices \eqref{anothermonodromy} we introduced
in this section as \\
$\displaystyle
{}_n \langle 0^k 1^{n-k}|
\prod_{j=k+1}^m \overline{A}(w_j|u_1,\dots,u_{n})
\prod_{j=1}^k \overline{C}(w_j|u_1,\dots,u_{n})
| \overline{\Omega} \rangle_n
$, and we have
\begin{align}
&W_{m+n-k,n}(u_1,\dots,u_n|w_1,\dots,w_{m+n-k}|1,\dots,k,m+1,\dots,m+n-k)
\nonumber \\
=&\prod_{j=1}^n u_j^{n-k}
{}_n \langle 0^k 1^{n-k}|
\prod_{j=k+1}^m \overline{A}(w_j|u_1,\dots,u_{n})
\prod_{j=1}^k \overline{C}(w_j|u_1,\dots,u_{n})
| \overline{\Omega} \rangle_n. \label{forrectangular}
\end{align}
Next, we commute the operators $\prod_{j=k+1}^m \overline{A}(w_j|u_1,\dots,u_{n})$
with
$\prod_{j=1}^k \overline{C}(w_j|u_1,\dots,u_{n})$.
The intertwining relation
\begin{align}
&R_{12}(w_1,w_2)
\overline{T}_2(w_2|u_1,\dots,u_n) \overline{T}_1(w_1|u_1,\dots,u_n)
\nonumber \\
=&
\overline{T}_1(w_1|u_1,\dots,u_n) \overline{T}_2(w_2|u_1,\dots,u_n) 
R_{12}(w_1,w_2),
\end{align}
gives
\begin{align}
\overline{A}(w_1) \overline{C}(w_2)
&=\frac{w_2}{w_2-w_1} \overline{C}(w_2) \overline{A}(w_1)
-\frac{w_1}{w_2-w_1} \overline{C}(w_1) \overline{A}(w_2), \\
\overline{A}(w_1) \overline{C}(w_2)
&=\overline{A}(w_2) \overline{C}(w_1), \\
\overline{A}(w_1) \overline{A}(w_2)
&=\overline{A}(w_2) \overline{A}(w_1), \\
\overline{C}(w_1) \overline{C}(w_2)
&=\overline{C}(w_2) \overline{C}(w_1),
\end{align}
from which we can get the following compact form of the commutation relation
by the argument in \cite{ShigechiUchiyama}

\begin{figure}[ht]
\includegraphics[width=12cm]{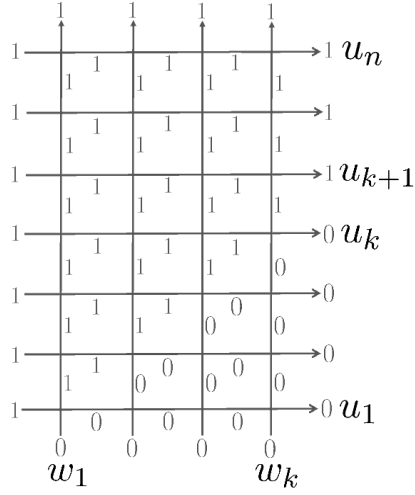}
\caption{
The partition function
$\displaystyle {}_n \langle 0^k 1^{n-k}|
\prod_{j \in S_k^m} \overline{C}(w_j|u_1,\dots,u_{n}) | \overline{\Omega} \rangle_n$.
We note that only one configuration is allowed,
from which one can see $\displaystyle {}_n \langle 0^k 1^{n-k}|
\prod_{j \in S_k^m} \overline{C}(w_j|u_1,\dots,u_{n}) | \overline{\Omega} \rangle_n=\prod_{j=1}^n u_j^k$.
}
\label{picturewavefunctionrrectangularfreeze}
\end{figure}

\begin{align}
&\prod_{j=k+1}^m \overline{A}(w_j|u_1,\dots,u_{n})
\prod_{j=1}^k \overline{C}(w_j|u_1,\dots,u_{n}) \nonumber \\
=&\sum_{S_k^m \in \binom{[m]}{k}}
\prod_{i \in S_k^m, j \in \overline{S_k^m}}
\frac{w_i}{w_i-w_j}
\prod_{j \in S_k^m} \overline{C}(w_j|u_1,\dots,u_{n})
\prod_{j \in \overline{S_k^m}} \overline{A}(w_j|u_1,\dots,u_{n}).
\label{multipletwo}
\end{align}
Using \eqref{multipletwo}
and the action of the $A$-operators on the state $| \overline{\Omega} \rangle_n$\begin{align}
\prod_{j \in \overline{S_k^m}} \overline{A}(w_j|u_1,\dots,u_{n})
| \overline{\Omega} \rangle_n
=\prod_{j \in \overline{S_k^m}} \prod_{i=1}^n (u_i-w_j)
| \overline{\Omega} \rangle_n,
\end{align}
\eqref{forrectangular} can be rewritten as
\begin{align}
&W_{m+n-k,n}(u_1,\dots,u_n|w_1,\dots,w_{m+n-k}|1,\dots,k,m+1,\dots,m+n-k)
\nonumber \\
=&\prod_{j=1}^n u_j^{n-k}
\sum_{S_k^m \in \binom{[m]}{k}}
\prod_{i \in S_k^m, j \in \overline{S_k^m}}
\frac{w_i}{w_i-w_j}
\prod_{j \in \overline{S_k^m}} \prod_{i=1}^n (u_i-w_j)
{}_n \langle 0^k 1^{n-k}|
\prod_{j \in S_k^m} \overline{C}(w_j|u_1,\dots,u_{n})
| \overline{\Omega} \rangle_n. \label{forrectangulartwo}
\end{align}

One can easily see from the graphical description that
the partition function \\
$\displaystyle {}_n \langle 0^k 1^{n-k}|
\prod_{j \in S_k^m} \overline{C}(w_j|u_1,\dots,u_{n}) | \overline{\Omega} \rangle_n$
is completely frozen
(Figure \ref{picturewavefunctionrrectangularfreeze}),
and find that its explicit form is given by
\begin{align}
{}_n \langle 0^k 1^{n-k}|
\prod_{j \in S_k^m} \overline{C}(w_j|u_1,\dots,u_{n}) 
| \overline{\Omega} \rangle_n =\prod_{j=1}^n u_j^k. \label{forsubstitution}
\end{align}
Substituting \eqref{forsubstitution} into \eqref{forrectangulartwo},
we get
\begin{align}
&W_{m+n-k,n}(u_1,\dots,u_n|w_1,\dots,w_{m+n-k}|1,\dots,k,m+1,\dots,m+n-k)
\nonumber \\
=&\prod_{j=1}^n u_j^{n}
\sum_{S_k^m \in \binom{[m]}{k}}
\prod_{i \in S_k^m, j \in \overline{S_k^m}}
\frac{w_i}{w_i-w_j}
\prod_{j \in \overline{S_k^m}} \prod_{i=1}^n (u_i-w_j).
\label{usethisforidentity}
\end{align}
Finally, we compare the two expressions \eqref{forcomparisontwo}
and \eqref{usethisforidentity} to get
\begin{align}
&\prod_{j=1}^n \frac{1}{(1-z_j)^{m+n-k}}
G_{((m-k)^{n-k},0^k)}(\{z_1,\dots,z_n \}|\bs{\alpha}) \nonumber \\
=&\prod_{j=1}^n u_j^{n}
\sum_{S_k^m \in \binom{[m]}{k}}
\prod_{i \in S_k^m, j \in \overline{S_k^m}}
\frac{w_i}{w_i-w_j}
\prod_{j \in \overline{S_k^m}} \prod_{i=1}^n (u_i-w_j),
\end{align}
which, after using the translation rule
$z_j=1-u_j^{-1} \ (j=1,\dots,n)$, $\alpha_j=1-w_j \ (j=1,\dots, m+n-k)$,
becomes the identity
\begin{align}
G_{((m-k)^{n-k},0^k)}(\bs{z}|\bs{\alpha})=\sum_{S_k^m \in \binom{[m]}{k}}
\frac{\displaystyle \prod_{i \in {S_k^m}} (1-\alpha_i)^{m-k} \prod_{j \in \overline{S_k^m}}
\prod_{i=1}^n (z_i \oplus \alpha_j)}
{\displaystyle \prod_{i \in {S_k^m}} \prod_{j \in \overline{{S_k^m}}}(\alpha_j-\alpha_i)}.
\end{align}

\end{proof}

{\bf Example of Theorem \ref{rectangulartheorem}} \\
When $n=2$, $k=1$, $m=2$, the right hand side of \eqref{rectangularidentity}
is
\begin{align}
&\frac{(1-\alpha_1)(z_1 \oplus \alpha_2)(z_2 \oplus \alpha_2)}
{\alpha_2-\alpha_1}
+\frac{(1-\alpha_2)(z_1 \oplus \alpha_1)(z_2 \oplus \alpha_1)}
{\alpha_1-\alpha_2} \nonumber \\
=&(\alpha_1 \alpha_2-\alpha_1-\alpha_2+1)(z_1+z_2-z_1 z_2)+\alpha_1+\alpha_2
-\alpha_1 \alpha_2,
\end{align}
which gives the left hand side
$G_{(1,0)}(\{z_1,z_2 \}|\bs{\alpha})$. \\

Combining \eqref{rectangularidentity}
with the Guo-Sun identity \eqref{guosun} for the case $\lambda_1=\cdots=\lambda_k=0$ and using $G_{(0,\dots,0)}(\bs{z}_{S_k^n}|\bs{\alpha})=1$,
we get the following duality formula.
\begin{theorem}
The following identity holds:
\begin{align}
&\sum_{S_k^n \in \binom{[n]}{k}}
\frac{\displaystyle \prod_{i \in {S_k^n}} (1- z_i)^{n-k} \prod_{j \in \overline{S_k^n}}
\prod_{i=1}^m
(z_j \oplus \alpha_i)}{\displaystyle \prod_{i \in {S_k^n}} \prod_{j \in \overline{{S_k^n}}}(z_j-z_i)}
\nonumber \\
=&\sum_{S_k^m \in \binom{[m]}{k}}
\frac{\displaystyle \prod_{i \in {S_k^m}} (1-\alpha_i)^{m-k} \prod_{j \in \overline{S_k^m}}
\prod_{i=1}^n (z_i \oplus \alpha_j)}
{\displaystyle \prod_{i \in {S_k^m}} \prod_{j \in \overline{{S_k^m}}}(\alpha_j-\alpha_i)}.
\end{align}
\end{theorem}

\section{A $q$-deformation}
In this section, we discuss a $q$-deformation of the
Guo-Sun identity.
We follow the same procedure of computation
in section 3 done for the five-vertex model.
Now we consider the $U_q(\widehat{sl_2})$ six-vertex model whose
$R$-matrix is given by \eqref{rmatrix}.
Recall that the correspondence between the
wavefunctions of the six-vertex model and the symmetric functions
\eqref{symmetricfunction} are given by
\eqref{correspondencesixvertex},
which applied to the one
$W_{m+n-k,n}(u_1,\dots,u_n|w_1,\dots,w_{m+n-k}|x_1,\dots,x_k,m+1,\dots,m+n-k)$ we deal with in this section, becomes
\begin{align}
&W_{m+n-k,n}(u_1,\dots,u_n|w_1,\dots,w_{m+n-k}|x_1,\dots,x_k,m+1,\dots,m+n-k)
\nonumber \\
=&
F_{m+n-k,n}(u_1,\dots,u_n|w_1,\dots,w_{m+n-k}|x_1,\dots,x_k,m+1,\dots,m+n-k).
\end{align}

\begin{figure}[ht]
\includegraphics[width=12cm]{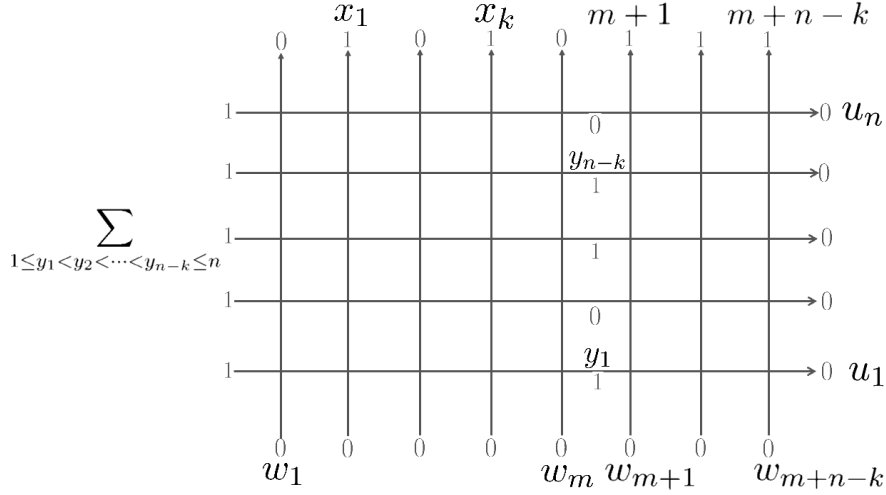}
\caption{The decomposition of the wavefunctions
\eqref{tochuqdeform}.
From the ice rule of the $R$-matrix of the six-vertex model
and the boundary condition, we note that for every allowed configuration,
the states on the horizontal edges
between the $m$-th and the $(m+1)$-th column consist of
$(n-k)$ 1's and $k$ 0's.
We label the positions of the $(n-k)$ 1's as $y_1, y_2, \dots, y_{n-k}$
($1 \le y_1 < y_2 < \cdots < y_{n-k} \le n$),
starting from the bottom row.
}
\label{pictureinsertcompleteness}
\end{figure}

Next we examine $W_{m+n-k,n}(u_1,\dots,u_n|w_1,\dots,w_{m+n-k}|x_1,\dots,x_k,m+1,\dots,m+n-k)$ from another point of view as in section 3.
Unlike the case for the five-vertex model,
there are several allowed configurations in the rightmost $n-k$ columns
(for the five-vertex model,
only one configuration is allowed).
However, from the
the so-called ice rule
$
{}_a \langle \gamma | {}_b \langle \delta | R_{a b}(u,w)
|\alpha \rangle_a | \beta \rangle_b=0$
unless $\alpha+\beta=\gamma+\delta$
and since the states at the right boundary are all $\langle 0|$,
one can see that the wavefunctions can be decomposed as
(Figure \ref{pictureinsertcompleteness})
\begin{align}
&W_{m+n-k,n}(u_1,\dots,u_n|w_1,\dots,w_{m+n-k}|x_1,\dots,x_k,m+1,\dots,m+n-k)
\nonumber \\
=&\sum_{1 \le y_1 < y_2 < \cdots < y_{n-k} \le n }
\overline{W}_{n,n-k}(w_{m+1},\dots,w_{m+n-k}|u_1,\dots,u_n|y_1,\dots,y_{n-k})
\nonumber \\
\times&
{}_{m} \langle x_1, \dots, x_k|
\prod_{j=y_{n-k}+1}^n B(u_j|w_1,\dots,w_m) D(u_{y_{n-k}}|w_1,\dots,w_m)
\prod_{j=y_{n-k-1}+1}^{y_{n-k}-1} B(u_j|w_1,\dots,w_m) \nonumber \\
\times& \cdots \times
\prod_{j=y_{1}+1}^{y_2-1} B(u_j|w_1,\dots,w_m) D(u_{y_1}|w_1,\dots,w_m)
\prod_{j=1}^{y_{1}-1} B(u_j|w_1,\dots,w_m)
| \Omega \rangle_m, \label{tochuqdeform}
\end{align}
where
\begin{align}
&\overline{W}_{n,n-k}(w_{m+1},\dots,w_{m+n-k}|u_1,\dots,u_n|y_1,\dots,y_{n-k})
\nonumber \\
=&{}_{m+1} \langle 1| \otimes \cdots \otimes {}_{m+n-k} \langle 1|
\prod_{j=y_{n-k}+1}^n A(u_j|w_{m+1},\dots,w_{m+n-k}) B(u_{y_{n-k}}|w_{m+1},\dots,w_{m+n-k})
\nonumber \\
\times&\prod_{j=y_{n-k-1}+1}^{y_{n-k}-1} A(u_j|w_{m+1},\dots,w_{m+n-k}) \cdots \times
\prod_{j=y_{1}+1}^{y_2-1} A(u_j|w_{m+1},\dots,w_{m+n-k}) \nonumber \\
\times&B(u_{y_1}|w_{m+1},\dots,w_{m+n-k})
\prod_{j=1}^{y_{1}-1} A(u_j|w_{m+1},\dots,w_{m+n-k})
| 0 \rangle_{m+n-k} \otimes \cdots \otimes |0 \rangle_{m+1}.
\end{align}
$\overline{W}_{n,n-k}(w_{m+1},\dots,w_{m+n-k}|u_1,\dots,u_n|y_1,\dots,y_{n-k})$
is another class of wavefunctions
and one can show the following correspondence
\begin{align}
&\overline{W}_{n,n-k}(w_{m+1},\dots,w_{m+n-k}|u_1,\dots,u_n|y_1,\dots,y_{n-k})
\nonumber \\
=&\overline{F}_{n,n-k}(w_{m+1},\dots,w_{m+n-k}
|u_1,\dots,u_n|y_1,\dots,y_{n-k}),
\label{anothercorrespondence}
\end{align}
where $\overline{F}_{n,n-k}(w_{m+1},\dots,w_{m+n-k}|u_1,\dots,u_n|y_1,\dots,y_{n-k})$ is the following symmetric functions
with symmetric variables $w_{m+1},\dots,w_{m+n-k}$,
another set of variables $u_1,\dots,u_n$
and a set of integers $y_1,\dots,y_{n-k}$
satisfying $1 \le y_1 < y_2 < \cdots < y_{n-k} \le n$
\begin{align}
&\overline{F}_{n,n-k}(w_{m+1},\dots,w_{m+n-k}|u_1,\dots,u_n|y_1,\dots,y_{n-k})
\nonumber \\
=&
\sum_{\sigma \in S_{n-k}}
\prod_{j=1}^{n-k} \prod_{i=y_j+1}^{n} q(u_{i}-w_{m+\sigma(j)})
\prod_{1 \le i < j \le n-k}
\frac{qw_{m+\sigma(i)}-w_{m+\sigma(j)}}{q(w_{m+\sigma(i)}-w_{m+\sigma(j)})}
\nonumber \\
&\times
\prod_{j=1}^{n-k} \prod_{i=1}^{y_j-1}(u_i-q w_{m+\sigma(j)})
\prod_{j=1}^{n-k} (1-q)u_{y_j}.
\label{anothersymmetricfunction}
\end{align}
One can show the correspondence \eqref{anothercorrespondence} for example
by the Izergin-Korepin method \cite{Ko,Iz},
which can be applied to the wavefunctions \cite{Motegi} as follows.
First, we construct the following Korepin's lemma
which list the properties of the partition functions
which uniquely characterize them.
For the case of the wavefunctions
$\overline{W}_{n,n-k}(w_{m+1},\dots,w_{m+n-k}
|u_1,\dots,u_n|y_1,\dots,y_{n-k})$,
the Korepin's Lemma is given below.
\begin{proposition} \label{korepinlemma}
The partition functions
$\overline{W}_{n,n-k}(w_{m+1},\dots,w_{m+n-k}|u_1,\dots,u_n|y_1,\dots,y_{n-k})$
satisfies the following properties. \\
\\
 (1) $\overline{W}_{n,n-k}(w_{m+1},\dots,w_{m+n-k}|u_1,\dots,u_n|y_1,\dots,y_{n-k})$ is a polynomial of degree $n-k$ in $u_n$
if $y_{n-k}=n$. \\
\\
 (2) $\overline{W}_{n,n-k}(w_{m+1},\dots,w_{m+n-k}|u_1,\dots,u_n|y_1,\dots,y_{n-k})$ is symmetric
with respect to $w_j$, $j=m+1,\dots,m+n-k$. \\
\\
(3) The following recursive relations between the
partition functions hold if $y_{n-k}=n$:
\begin{align}
&\overline{W}_{n,n-k}(w_{m+1},\dots,w_{m+n-k}|u_1,\dots,u_n|y_1,\dots,y_{n-k})|_{u_n=0}=0, \\
&\overline{W}_{n,n-k}(w_{m+1},\dots,w_{m+n-k}|u_1,\dots,u_n|y_1,\dots,y_{n-k})|_{u_n=w_{m+n-k}}
\nonumber \\
=&(1-q) w_{m+n-k} \prod_{j=m+1}^{m+n-k-1} (w_{m+n-k}-q w_j)
\prod_{j=1}^{n-1}(u_j-q w_{m+n-k})
\nonumber \\
&\times \overline{W}_{n-1,n-k-1}(w_{m+1},\dots,w_{m+n-k-1}
|u_1,\dots,u_{n-1}|y_1,\dots,y_{n-k-1}).
\end{align}
If $y_{n-k} \neq n$, the following factorizations hold for the wavefunctions:
\begin{align}
&\overline{W}_{n,n-k}(w_{m+1},\dots,w_{m+n-k}|u_1,\dots,u_n|y_1,\dots,y_{n-k})
\nonumber \\
=&\prod_{j=m+1}^{m+n-k} q(u_n-w_j)
\overline{W}_{n-1,n-k}(w_{m+1},\dots,w_{m+n-k}|u_1,\dots,u_{n-1}|y_1,\dots,y_{n-k}).
\end{align}
\\
(4) The following holds for the case $n-k=1$, $y_1=n$
\begin{align}
&
\overline{W}_{n,1}(w_{m+1}|u_1,\dots,u_n|n)
=(1-q)u_n \prod_{j=1}^{n-1} (u_j-q w_{m+1}).
\end{align}
\end{proposition}
After constructing Korepin's Lemma, by showing
that the symmetric functions \\
$\overline{F}_{n,n-k}(w_{m+1},\dots,w_{m+n-k}|u_1,\dots,u_n|y_1,\dots,y_{n-k})$
satisfy all the properties in Proposition \ref{korepinlemma}
and we get the correspondence \eqref{anothercorrespondence}.

Now we examine the factors
\begin{align}
&
{}_{m} \langle x_1, \dots, x_k|
\prod_{j=y_{n-k}+1}^n B(u_j|w_1,\dots,w_m) D(u_{y_{n-k}}|w_1,\dots,w_m)
\prod_{j=y_{n-k-1}+1}^{y_{n-k}-1} B(u_j|w_1,\dots,w_m) \nonumber \\
\times& \cdots \times
\prod_{j=y_{1}+1}^{y_2-1} B(u_j|w_1,\dots,w_m) D(u_{y_1}|w_1,\dots,w_m)
\prod_{j=1}^{y_{1}-1} B(u_j|w_1,\dots,w_m)
| \Omega \rangle_m, \label{beforesimplifying}
\end{align}
in \eqref{tochuqdeform}.
We apply the method used in
\cite{KMT} to study correlation functions of the
XXZ spin chain for simplifying
\eqref{beforesimplifying}.
From the intertwining relation \eqref{rttone}
for the six-vertex model \eqref{rmatrix},
we get
\begin{align}
D(u_1)B(u_2)&=\frac{u_1-q u_2}{u_1-u_2}B(u_2)D(u_1)
+\frac{(q-1)u_2}{u_1-u_2}B(u_1)D(u_2), \label{qrttone} \\
B(u_1)B(u_2)&=B(u_2)B(u_1). \label{qrtttwo}
\end{align}
From the argument which is standard in the algebraic Bethe ansatz,
we can combine \eqref{qrttone}, \eqref{qrtttwo} and
\begin{align}
D(u|w_1,\dots,w_m)|\Omega \rangle_m
=\prod_{i=1}^m (u-w_i)|\Omega \rangle_m,
\end{align}
to show the following relation
\begin{align}
&D(u_{\ell+1}|w_1,\dots,w_m) \prod_{j=1}^\ell B(u_j|w_1,\dots,w_m) | \Omega \rangle_m \nonumber \\
=&\sum_{k=1}^{\ell+1} \prod_{i=1}^m (u_k-w_i)
\frac{\displaystyle \prod_{j=1}^\ell (u_k-q u_j)}{\displaystyle
\prod_{\substack{j=1 \\ j \neq k}}^{\ell+1} (u_k-u_j)}
\prod_{\substack{j=1 \\ j \neq k}}^{\ell+1} B(u_j|w_1,\dots,w_m) | \Omega \rangle_m
\label{BonD}.
\end{align}
Using \eqref{BonD} repeatedly, one gets
\begin{align}
&
\prod_{j=y_{n-k}+1}^n B(u_j|w_1,\dots,w_m) D(u_{y_{n-k}}|w_1,\dots,w_m)
\prod_{j=y_{n-k-1}+1}^{y_{n-k}-1} B(u_j|w_1,\dots,w_m) \nonumber \\
\times& \cdots \times
\prod_{j=y_{1}+1}^{y_2-1} B(u_j|w_1,\dots,w_m) D(u_{y_1}|w_1,\dots,w_m)
\prod_{j=1}^{y_{1}-1} B(u_j|w_1,\dots,w_m)
| \Omega \rangle_m \nonumber \\
=&\sum_{a_1=1}^{y_1} \sum_{\substack{a_2=1 \\ a_2 \neq a_1}}^{y_2}
\cdots \sum_{\substack{a_{n-k}=1 \\ a_{n-k} \neq a_1,\dots,a_{n-k-1}}}
^{y_{n-k}} \prod_{j=1}^{n-k}
\prod_{i=1}^m (u_{a_j}-w_i)
\frac{\displaystyle \prod_{j=1}^{n-k} \prod_{\substack{b=1 \\ b \neq a_1,\dots,a_{j-1}}}^{y_j-1} (u_{a_j}-q u_b)}{\displaystyle \prod_{j=1}^{n-k} \prod_{\substack{b=1 \\ b \neq a_1,\dots,a_j}}
^{y_j} (u_{a_j}-u_b)}
\nonumber \\
\times&\prod_{\substack{j=1 \\ j \neq a_1,\dots,a_{n-k}}}^{n} B(u_j|w_1,\dots,w_m) | \Omega \rangle_m. \label{qmultipleaction}
\end{align}
Combining \eqref{qmultipleaction} and the correspondence
\begin{align}
&{}_m \langle x_1,\dots,x_k|
\prod_{\substack{j=1 \\ j \neq a_1,\dots,a_{n-k}}}^{n} B(u_j|w_1,\dots,w_m) | \Omega \rangle_m \nonumber \\
=&W_{m,n}(\{u_1,\dots,u_n \} \backslash \{ u_{a_1}, \dots,u_{a_{n-k}} \}|w_1,\dots,w_m|x_1,\dots,x_k),
\end{align}
\eqref{beforesimplifying} becomes
\begin{align}
&
{}_{m} \langle x_1, \dots, x_k|
\prod_{j=y_{n-k}+1}^n B(u_j|w_1,\dots,w_m) D(u_{y_{n-k}}|w_1,\dots,w_m)
\prod_{j=y_{n-k-1}+1}^{y_{n-k}-1} B(u_j|w_1,\dots,w_m) \nonumber \\
\times& \cdots \times
\prod_{j=y_{1}+1}^{y_2-1} B(u_j|w_1,\dots,w_m) D(u_{y_1}|w_1,\dots,w_m)
\prod_{j=1}^{y_{1}-1} B(u_j|w_1,\dots,w_m)
| \Omega \rangle_m \nonumber \\
=&
\sum_{a_1=1}^{y_1} \sum_{\substack{a_2=1 \\ a_2 \neq a_1}}^{y_2}
\cdots \sum_{\substack{a_{n-k}=1 \\ a_{n-k} \neq a_1,\dots,a_{n-k-1}}}
^{y_{n-k}} \prod_{j=1}^{n-k}
\prod_{i=1}^m (u_{a_j}-w_i)
\frac{\displaystyle \prod_{j=1}^{n-k} \prod_{\substack{b=1 \\ b \neq a_1,\dots,a_{j-1}}}^{y_j-1} (u_{a_j}-q u_b)}{\displaystyle \prod_{j=1}^{n-k} \prod_{\substack{b=1 \\ b \neq a_1,\dots,a_j}}
^{y_j} (u_{a_j}-u_b)}
\nonumber \\
\times&W_{m,n}(\{u_1,\dots,u_n \} \backslash \{ u_{a_1}, \dots,u_{a_{n-k}} \}|w_1,\dots,w_m|x_1,\dots,x_k).
\label{aftersimplifying}
\end{align}
Inserting \eqref{aftersimplifying} into the right hand side of
\eqref{tochuqdeform}, one gets
\begin{align}
&W_{m+n-k,n}(u_1,\dots,u_n|w_1,\dots,w_{m+n-k}|x_1,\dots,x_k,m+1,\dots,m+n-k)
\nonumber \\
=&\sum_{1 \le y_1 < y_2 < \cdots < y_{n-k} \le n }
\sum_{a_1=1}^{y_1} \sum_{\substack{a_2=1 \\ a_2 \neq a_1}}^{y_2}
\cdots \sum_{\substack{a_{n-k}=1 \\ a_{n-k} \neq a_1,\dots,a_{n-k-1}}}
^{y_{n-k}} \prod_{j=1}^{n-k}
\prod_{i=1}^m (u_{a_j}-w_i)
\frac{\displaystyle \prod_{j=1}^{n-k} \prod_{\substack{b=1 \\ b \neq a_1,\dots,a_{j-1}}}^{y_j-1} (u_{a_j}-q u_b)}{\displaystyle \prod_{j=1}^{n-k} \prod_{\substack{b=1 \\ b \neq a_1,\dots,a_j}}
^{y_j} (u_{a_j}-u_b)} \nonumber \\
\times&\overline{W}_{n,n-k}(w_{m+1},\dots,w_{m+n-k}|u_1,\dots,u_n|y_1,\dots,y_{n-k}) \nonumber \\
\times&W_{m,n}(\{u_1,\dots,u_n \} \backslash \{ u_{a_1}, \dots,u_{a_{n-k}} \}|w_1,\dots,w_m|x_1,\dots,x_k),
\end{align}
which, using the correspondences
\eqref{correspondencesixvertex} and
\eqref{anothercorrespondence}, becomes
an identity
\begin{align}
&F_{m+n-k,n}(u_1,\dots,u_n|w_1,\dots,w_{m+n-k}|x_1,\dots,x_k,m+1,\dots,m+n-k)
\nonumber \\
=&\sum_{1 \le y_1 < y_2 < \cdots < y_{n-k} \le n }
\sum_{a_1=1}^{y_1} \sum_{\substack{a_2=1 \\ a_2 \neq a_1}}^{y_2}
\cdots \sum_{\substack{a_{n-k}=1 \\ a_{n-k} \neq a_1,\dots,a_{n-k-1}}}
^{y_{n-k}} \prod_{j=1}^{n-k}
\prod_{i=1}^m (u_{a_j}-w_i)
\frac{\displaystyle \prod_{j=1}^{n-k} \prod_{\substack{b=1 \\ b \neq a_1,\dots,a_{j-1}}}^{y_j-1} (u_{a_j}-q u_b)}{\displaystyle \prod_{j=1}^{n-k} \prod_{\substack{b=1 \\ b \neq a_1,\dots,a_j}}
^{y_j} (u_{a_j}-u_b)} \nonumber \\
\times&\overline{F}_{n,n-k}(w_{m+1},\dots,w_{m+n-k}|u_1,\dots,u_n|y_1,\dots,y_{n-k}) \nonumber \\
\times&F_{m,n}(\{u_1,\dots,u_n \} \backslash \{ u_{a_1}, \dots,u_{a_{n-k}} \}|w_1,\dots,w_m|x_1,\dots,x_k),
\end{align}
for the symmetric functions \eqref{symmetricfunction}
and \eqref{anothersymmetricfunction}.

\section{Conclusion}
In this paper, from the point of view of quantum integrability,
we first investigated the identity
for the factorial Grothendieck polynomials found by Guo and Sun \cite{GS}
which generalizes the one for the Schur polynomials by
Feh\'er, N\'emethi and Rim\'anyi \cite{FNR}.
We gave another proof by using the quantum inverse scattering method,
which is a method to analyze quantum integrable models.
Why the method can be used is based on the fact
between the correspondence between the wavefunctions
of a five-vertex model and the factorial Grothendieck polynomials.
See \cite{MS,MS2,GK,WZ} also for previous works on the investigations
of Cauchy-type identities, Gromov-Witten invariants and the
Littlewood-Richardson coefficients using this correspondence.

We next used the same idea and technique 
``in another direction'' to derive an identity
for the factorial Grothendieck poylnomials of rectangular shapes.
Combining the identity with the Guo-Sun identity, we
obtained a duality formula.
We also discussed a $q$-deformation of the Guo-Sun identity,
based on the correspondence between the wavefunctions
of the $U_q(\widehat{sl_2})$ six-vertex model and the $q$-deformation of the factorial
Grothendieck polynomials
and following the same line of computation to prove the Guo-Sun identity.
The identity obtained for the $q$-deformed symmetric functions
is rather much more complicated than the
Guo-Sun identity since the six-vertex model
is a more general model than the five-vertex model,
and it is an interesting
problem whether one can simplify the identity
to a more compact form.

It may also be interesting to reexamine existing formulas
for the Schur and Grothendieck polynomials from the
viewpoint of quantum integrability,
and apply the same idea and technique in different ways,
cases and places to obtain new identities.
It is also interesting to investigate if the integrability
technique can be applied beyond
Grassmannian Grothendieck polynomials.
One needs first to investigate if
the set-valued tableaux descriptions of
more general Grothendieck polynomials in
\cite{Matsumura,MatSugi,FanGuo} can be
translated into the language of integrable models.

\section*{Acknowledgments}
The author thanks the referee for careful reading
the manuscript and useful comments and suggestions.
This work was partially supported by grant-in-Aid
for and Scientific Research (C) No. 18K03205 and No. 16K05468.

\bibliographystyle{plainnat}

\end{document}